# A Human Story of Curiosity and Relevance

**Claire Burrin (University of Zurich)**

Mathematics is changing. Computers are verifying proofs, checking calculations, and exploring complex structures that would overwhelm human effort. Nonetheless, human creativity remains irreplaceable: we ask the questions and imagine the possibilities. This partnership is transforming research and also has consequences for teaching.

Computers can perform calculations and check proofs, but the responsibility to understand, question, validate, is still ours. Mathematical thinking helps to distinguish reliable knowledge from mere outputs. This discipline of rigour is indispensable to ensure that our increasing reliance on black boxes (computers, algorithms, AI) enhances rather than diminishes our intellectual judgment. And so, we will be increasingly tasked with guiding students in interpreting outputs, asking the meaningful questions, and appreciating the conceptual and creative side of mathematics.

Why do mathematicians study some problems and not others? What were the original questions that inspired their work? How do abstract ideas, sometimes centuries later, shape the world we live in? Here, digital security, data transmission, compression, mathematical modelling, all come to mind.

Of course, other sciences share this trajectory. A physicist might point to quantum mechanics, once a purely theoretical exercise, now central to semiconductors and lasers. Or a chemist might point to the periodic table, now foundational for drug design and materials. Curiosity-driven research is where tomorrow's breakthroughs are quietly prepared. Making this point clear to students — the next generation — is also part of the mission of a strong mathematical education.

I want to illustrate these ideas with a few anecdotes about a mathematical object that is fundamental to my own research, the lattice.

**From Babylon to Digital Security**

Whenever we send an email or use our credit card online, we rely on cryptographic protocols to keep our data secure. One of the pillars of internet security is the RSA algorithm, and its principle is taught in every algebra course: every whole number can be uniquely factored as a product of prime numbers — the so-called fundamental theorem of arithmetic.

The fascination with factoring is ancient. Clay tablets from ancient Babylonia, dating back to nearly four thousand years, record sophisticated tables of reciprocals and products, hinting at an early understanding of how numbers can be factored and combined. Aristotle reflected on those numbers that cannot be divided further, the prime numbers. Primes are, in effect, the atoms of numbers, simple building blocks from which all whole numbers are constructed.

Multiplying two large primes is easy, a computer does it in milliseconds. However, reversing the process, to factor back the product into primes, is a completely different story: once the primes reach a few hundred digits, the best-known algorithms take exponential time to solve this problem.

Modern security rests on the assumption that such a problem will remain hard for computers. Yet we know since the 1990s that a quantum computer could, in principle, factor whole numbers much more efficiently, rendering RSA and similar protocols obsolete. If building a quantum computer remains a great technological challenge, the uncertainty surrounding it has accelerated the race for post-quantum cryptography. Among the most promising candidates, we find lattice-based protocols.

In simple terms, a lattice is a regular grid of points extending through space. The definition can be explained to a student of linear algebra. Pick a set

Claire Burrin (University of Zurich)

of linearly independent vectors $v_1, v_2, \ldots, v_n$ in $\mathbb{R}^n$ and take all possible integer linear combinations

$$\mathbb{Z}v_1 + \mathbb{Z}v_2 + \ldots + \mathbb{Z}v_n$$

of these vectors. This produces a lattice: an orderly structure that exists in any number of dimensions.

The simplest and most famous hard problem associated with lattice-based cryptography is the Shortest Vector Problem (SVP): given a lattice, find its (nonzero) vector of shortest length. In two or three dimensions, it is easy: you just look. Once in thousands or millions of dimensions, the base vectors can point in many strange directions, making it exceptionally difficult to predict which integer combination results in the shortest vector. Finding it amounts to locating the proverbial needle in a multidimensional haystack.

Where a basic concept like prime numbers was enough for the first rudimentary encryption protocols, high-level modern number theory is at play to secure the most powerful protocols. For example, at a workshop I was co-organizing at the SwissMAP Research Station in Les Diablerets[1] this summer, a new robustness result for SVP was announced in connection to deep questions in number theory such as the Grand Riemann Hypothesis.

**The Curse of Dimensionality**

There's a joke among scientists: the engineer asks the mathematician, "How can you picture eleven-dimensional space?" The mathematician shrugs. "I imagine n-dimensional space and let n equal eleven." How often has mathematics been accused of being esoteric! But in applications to the Information Age, working in very large dimensions is actually the norm. A machine learning algorithm processes data as vectors in spaces with millions of coordinates. Search engines, when they rank pages, manipulate gigantic matrices. In fact, the search engine Google got its name from the googol[2] $10^{100}$ which is meant to represent an unimaginably large number, far greater than the estimated number of atoms in the observable universe.

It then makes sense to develop a picture for eleven-dimensional, or one billion-dimensional space. In fact, it turns out that higher dimensional geometry behaves in ways that defy our everyday 2D/3D intuition. For example, in a high-dimensional cube, most of the volume concentrates not at the centre but rather near the corners, which is the opposite of what we observe in low dimensions. Many problems also see an exponential growth in complexity as dimension increases. To see this, imagine trying to sample the unit square with evenly spaced points — lattice points — at a distance of 0.1 of each other. For this we need one hundred points, but to do the same in 9 dimensions will require one billion points.

Another striking example is the sphere packing problem: What is the most space-efficient way to arrange spheres of equal size in space?[3] In three dimensions, the best possible arrangement is the one you would expect: cannonball packing, that is, the pyramidal stacking of cannonballs. But in (most) higher dimensions, we do not even have an idea of what the densest packing would look like. Here again, we have an abstract question with very concrete consequences; in this case to coding theory. Imagine transmitting a message across a noisy channel. Each message can be represented as a point in high-dimensional space, and noise slightly shifts the point. By surrounding each point with a sphere representing the range of noise that can be tolerated, the received signal can still be decoded correctly. The goal is then to pack these spheres as densely as possible so that more information can be reliably transmitted.

Despite rapid and exciting developments in the last few years, still very little is known about the sphere packing problem in higher dimensions. Earlier this year, Boaz Klartag of the Weizmann Institute — currently visiting ETH for the Fall semester — announced a remarkable result based on combining geometry and probability theory: the existence of a high-dimensional lattice sphere packing that is far denser than that was known before. Can we exhibit a concrete lattice with this property? Could these constructions improve coding protocols or other technologies?

---

[1] Cf. https://swissmaprs.ch/.

[2] The term *googol*—essentially a made-up name—has a somewhat debated origin, as explained in sources such as Wikipedia.org/wiki/Googol, which connects interestingly to the theme of popular mathematics.

[3] Cf. Einstein Lectures: Maryna Viazovska - The sphere packing problem.

**Proof and Certainty**

Abstract thinking is not unique to mathematics; physicists, engineers, or computer scientists all develop their own forms of conceptual reasoning. Mathematics is different in one essential respect: its currency is not experiment, but proof.

A proof is a rigorous chain of reasoning that shows a statement must be true. Unlike in other sciences, a proven theorem does not get overturned; its truth is permanent. This tradition goes back to the ancient Greeks, taking shape in Euclid's Elements (300 BC) and continues to define mathematical practice today. So much so that a mathematician's career is measured by what kind of theorems they prove, how many, and the elegance of their reasoning.

The story of the resolution of Kepler's conjecture, that the cannonball packing is the most efficient arrangement of spheres in 3D, is a very good illustration of the changes under way in the practice of mathematical research. If the balls have to be organised along a lattice, then Kepler's conjecture is easy to prove. But could a small perturbation of such a rigid arrangement produce yet a more space-efficient packing? Centuries of intuition could not settle the question.

It is the advent of the computer in the 20[th] century that finally did so. László Féjes Tóth observed in the 1950s that Kepler's conjecture could in principle be reduced to a finite, but exceedingly large, number of calculations. In the early 1990s, Thomas Hales and his student Samuel Ferguson carried out this program: a computer-based proof of the Kepler conjecture by case-exhaustion[4]. The result was a manuscript shy of 300 pages, complemented by 40000 lines of codes, solving about 100000 linear programming problems.

At the time, assessing correctness meant painstakingly checking each detail by hand. A jury of twelve experts spent five years reviewing the work and ultimately declared to be "99% certain" of its correctness… a verdict that would have felt frustrating to any mathematician. Two decades later, the original paper was complemented with a formal, that is computer-verifiable, proof.

The nature of proof might be evolving but mathematics remains what it has always been: a discipline of curiosity. Every theorem begins as a question, sometimes naïve or purely theoretical, posed by someone who wanted to understand the symmetry of a form, or the shape of a pattern. Curiosity sparks understanding, which in turn enables applications — some that we cannot yet imagine.

University of Zurich, Institute of Mathematics, Winterthurerstrasse 190, CH-8057 Zurich.

E-mail: claire.burrin@math.uzh.ch

Website: https://user.math.uzh.ch/burrin/

ORCID: 0000-0003-0036-5539

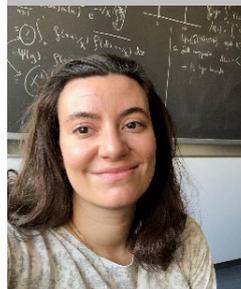

**Claire Burrin**, Dr. sc. ETH, is assistant professor at the Institute of Mathematics of the University of Zurich. Her research lies at the interface of number theory and dynamical systems, motivated by questions such as: How do arithmetic and geometric structures give rise to both random-like and extremal phenomena, and which analytic tools can be used to describe them?

She was previously a Hill Assistant Professor at Rutgers University in New Jersey, a senior assistant at ETH Zurich, and held various visiting positions, including in Princeton, Berkeley, and Bonn. She is part of the Zurich Dynamics Group, the NCCR SwissMAP—The Mathematics of Physics, and the Swiss-French collaborative project Equidistribution in Number Theory.

*Photo credentials: Private*

---

[4] References to Hales's proof of the Kepler conjecture include the full scientific article in the *Annals of Mathematics* (https://annals.math.princeton.edu/wp-content/uploads/annals-v162-n3-p01.pdf) and a more accessible overview in his *Notices of the AMS* article (https://www.ams.org/notices/200004/fea-hales.pdf).